\documentclass{article}

\usepackage{amssymb}
\usepackage{graphicx}
\begin{document}

\title{A naive procedure for computing angular spheroidal functions}

\author{J Sesma\footnote{javier@unizar.es}\\  \   \\
Departamento de F\'{\i}sica Te\'orica, Facultad de Ciencias,
\\ 50009 Zaragoza, Spain. \\  \ }

\maketitle

%\keywords{angular spheroidal wave equation; eigenvalues; eigenfunctions}
MSC [2010] {  33E10 (Primary)   33F05, 34L16, 65D20 (Secondary)}

PACS{ 02.30.Hq \and 02.30.Gp \and 03.65.Ge}

\begin{abstract}
An algorithm for computing eigenvalues and eigenfunctions of the angular spheroidal wave equation, based on a known but scarcely used method, is developed. By requiring the regularity of the wave function, represented by its series expansion, the eigenvalues appear as the zeros of a one variable function easily computable. The iterative extended Newton method is suggested as especially suitable for determining those zeros. The computation of the eigenfunctions is then immediate. The usefulness of the method, applicable also in the case of complex values of the ``prolateness" parameter, is illustrated by comparing its results with those of procedures used by other authors.
\end{abstract}

\section{Introduction}

The usefulness of spheroidal functions in many branches of Physics, like Quantum Mechanics, General Relativity, Signal Processing, etc., is well known and it does not need to be stressed. Due to that usefulness, the description of the spheroidal equation and of the main properties of its solutions deserves a chapter in handbooks of special functions like that by Abramowitz and Stegun \cite[Chap.~21]{abra}, the best known one, or the {\it NIST Digital Library of Mathematical Functions} \cite[Chap.~30]{nist}, the most recent one.

A review of the procedures used in the past century for obtaining the eigenvalues and eigenfunctions of the spheroidal wave equation can be found in a paper by Li {\em et al.} \cite{lilw}, where also an algorithm, implemented with the software package \verb"Mathematica", is provided.
In the present century, articles dealing with the solutions of the angular spheroidal wave equation have continued appearing. Without aiming to be exhaustive, let us mention the papers by Aquino {\em et al.}~\cite{aqui}, Falloon {\em et al.}~\cite{fall}, Boyd \cite{boyd}, Barrowes {\em et al.} \cite{barr}, Walter and Soleski \cite{walt}, Abramov and Kurochkin \cite{aaab}, Kirby \cite{kirb}, Karoui and Moumni \cite{karo}, Gosse \cite{goss}, Tian \cite{tian}, Rokhlin and Xiao \cite{rokh}, Osipov and Rokhlin \cite{osi1},  Ogburn {\em et al.}~\cite{ogbu} and Huang {\em et al.}~\cite{huan}, and the books by Hogan and Lakey \cite{hoga}, and by Osipov, Rokhlin and Xiao \cite{osi2}.

Different strategies have been used to solve the angular spheroidal wave equation. The classical procedure starts with the angular spheroidal wave function written as a series of solutions of another similar differential equation, commonly the Legendre one, with coefficients obeying a three term recurrence relation. The resulting expansion becomes convergent only when such coefficients constitute a minimal solution of the recurrence relation. The eigenvalue problem encountered in this way is solved either as a transcendental equation involving a continued fraction, or written in a matrix form. Procedures based on the direct solution of the angular spheroidal equation, without having recourse to comparison with other differential equations, have been less frequently used. The relaxation method proposed by Caldwell \cite{cald} and reproduced, as a worked example, in the  {\it Numerical Recipes} \cite[Sec. 17.4]{pres}, and the finite difference algorithm, described in the recently appeared paper by Ogburn {\em et al.} \cite{ogbu}, deserve to be quoted. Here we suggest to follow a procedure, based also on the direct treatment of the spheroidal equation, which benefits from an idea that can be found in a paper by Skorokhodov and Khristoforov \cite{skor} dealing with the singularities of the eigenvalues $\lambda_{m,n}$ considered as function of the (complex) prolateness parameter $c$. A shooting method is used. But, instead of imposing the boundary conditions to a numerically obtained solution, algebraic regular solutions around the regular point $\eta=0$ or around the regular singular point $\eta=1$ are written. Smooth matching of both solutions, i. e. cancelation of their Wronskian, at any point $\eta\in(-1, 1)$ determines the eigenvalues. In our implementation of the procedure, we choose $\eta=0$ as matching point.

A discomfort, when dealing with spheroidal wave functions, is the lack of universality of the notation used to represent them. The {\em Digital Library of Mathematical Functions} \cite[Chap.~30]{nist} provides information about the different notations found in the bibliography. Here we adopt, for the eigenvalues and eigenfunctions, the notation of the {\em Handbook of Mathematical Functions} \cite[Chap.~21]{abra}. The same notation is used in Ref. \cite{lilw}, a paper whose results we will try to reproduce, for comparison, with the method here developed.

In the next section, we recall the angular spheroidal equation and write its solutions in the form of power series expansions around the origin and around the singular point $\eta=1$. The procedure for computing the eigenvalues is presented in Section 3. The results of its application in some typical cases are also given. Section 4 shows that normalized eigenfunctions can be trivially obtained. Some figures illustrate the procedure. A few final comments are contained in Section 5.

\section{The differential equation}

The angular spheroidal wave function $S_{m,n}(c,\eta)$, defined in the interval $-1\leq \eta \leq 1$, satisfies the differential equation \cite[Eq. 21.6.2]{abra}
\begin{equation}
\frac{d}{d\eta}\left[(1-\eta^2)\frac{d}{d\eta}S_{m,n}(c,\eta)\right]+\left(\lambda_{m,n}-c^2\eta^2-\frac{m^2}{1-\eta^2}\right)S_{m,n}(c,\eta)=0 \label{ii1}
\end{equation}
stemming from the separation of the wave equation in spheroidal coordinates, with separation constants $m$ and  $\lambda_{m,n}$. Periodicity of the azimuthal part of the wave restricts the values of $m$ to the integers and, given the invariance of the differential equation in the reflection $m\Rightarrow -m$, only non-negative integer values of $m$ need to be considered. The other separation constant, $\lambda_{m,n}$, commonly referred to as eigenvalue, must be such that $S_{m,n}(c,\eta)$ becomes finite at the singular points $\eta=\pm 1$. Their different values, for given $m$ and $c^2$, are labeled  by the integer $n$. In most applications, the external parameter $c^2$ is real, positive in the case of prolate coordinates and negative for oblate ones. There are, however, interesting cases corresponding to complex values of $c^2$ \cite{aaab,barr,lilw,ogbu,oguc,skor}.

Instead of solving directly Eq. (\ref{ii1}), it is convenient to introduce the change of function
\begin{equation}
S_{m,n}(c,\eta)=(1-\eta^2)^{m/2}\,w(\eta)\,,  \label{ii2}
\end{equation}
and to solve the differential equation
\begin{equation}
(1-\eta^2)\,\frac{d^2}{d\eta^2}w(\eta)-2(m+1)\,\eta\,\frac{d}{d\eta}w(\eta) +\left(z-c^2\,\eta^2\right)w(\eta)=0\,, \label{ii3}
\end{equation}
where
\begin{equation}
z\equiv \lambda_{m,n}-m(m+1)  \label{ii4}
\end{equation}
is considered as the new eigenvalue.

Two independent solutions about the ordinary point $\eta=0$, valid in the interval $-1<\eta<1$, are
\begin{equation}
w_\sigma(\eta)=\sum_{k=0}^\infty\,a_{k,\sigma}\,\eta^{k+\sigma}\,, \qquad \sigma=0, 1\,, \label{ii5}
\end{equation}
with coefficients given by the recurrence relation
\begin{eqnarray}
a_{0,\sigma}=1\,,\quad a_{1,0}=0\,,\quad (k+\sigma)(k-1+\sigma)\,a_{k,\sigma}=& & \nonumber \\
& & \hspace{-160pt} \left[(k-1+2m+\sigma)(k-2+\sigma)-z\right]\,a_{k-2,\sigma}+c^2\,a_{k-4,\sigma}\,. \label{ii6}
\end{eqnarray}
Obviously, $w_0(\eta)$ and $w_1(\eta)$ are respectively  even and odd functions of $\eta$.

Solutions about the regular singular point $\eta=1$ can also be written. In terms of the variables
\begin{equation}
t\equiv 1-\eta\,,  \qquad u(t)\equiv w(1-\eta)\,,  \label{ii7}
\end{equation}
the differential equation (\ref{ii3}) turns into
\begin{equation}
t(2-t)\frac{d^2}{dt^2}u(t)+2(m+1)(1-t)\,\frac{d}{dt}u(t)+\left[z-c^2(1-t)^2\right]u(t)=0\,.  \label{ii8}
\end{equation}
The solution of this equation which makes $S_{m,n}$ to be regular at $t=0$ is,
except for an arbitrary multiplicative constant,
\begin{equation}
u_{\rm reg}(t)=\sum_{j=0}^\infty b_j\,t^j\,,   \label{ii9}
\end{equation}
with coefficients given by
\begin{eqnarray}
b_0=1\,,\qquad 2j(j+m)\,b_j=\left[(j-1)(j+2m)-z+c^2\right]b_{j-1}& &  \nonumber\\
& & \hspace{-80pt} -\,2c^2\,b_{j-2}+c^2\,b_{j-3}\,.   \label{ii10}
\end{eqnarray}
In terms of the variable $\eta$, this regular solution, valid for $-1<\eta\leq 1$, is
\begin{equation}
w_{\rm reg}(\eta)=u_{\rm reg}(1-\eta)\,.  \label{ii11}
\end{equation}

\section{The eigenvalues}

The problem of finding the eigenvalues $\lambda_{m,n}$, for given $c^2$ and $m$, reduces to require the regularity of $w_\sigma(\eta)$ at $\eta=1$. (The regularity at $\eta=-1$ is then implied by the symmetry of $w_\sigma$.) In the particular case of being $c^2=0$, the problem can be solved algebraically. The recurrence relation in (\ref{ii6}) reduces in this case to
\begin{equation}
(k+\sigma)(k-1+\sigma)\,a_{k,\sigma}=
 \left[(k-1+2m+\sigma)(k-2+\sigma)-z\right]\,a_{k-2,\sigma}\,. \label{iii1}\\
\end{equation}
Obviously, the series in the right hand side of (\ref{ii5}) is divergent for $\eta=1$ unless the value of $z$ is such that one of the $a_{k,\sigma}$ of even subindex, say $a_{2K+2,\sigma}$ (with $K=0, 1, 2, \ldots$), becomes zero, in which case $w_\sigma(\eta)$ turns out to be a polynomial of degree $2K+\sigma$. This happens for
\begin{equation}
z=(2K+1+2m+\sigma)(2K+\sigma). \label{iii2}
\end{equation}
By using the habitual notation
\begin{equation}
2K+\sigma \equiv n-m   \label{iii3}
\end{equation}
for the degree of the polynomial, we obtain for the eigenvalues in the case of $c^2=0$
\begin{equation}
z=(n+m+1)(n-m)\,,  \label{iii4}
\end{equation}
that is, in view of (\ref{ii4}),
\begin{equation}
\lambda_{m,n}(c^2=0)=n(n+1)\,.  \label{iii5}
\end{equation}

For $c^2\neq 0$, a convenient way of guaranteeing the regularity of $w_\sigma(\eta)$ at $\eta=1$ is to require the cancelation of the Wronskian $W$ of $w_\sigma$ and $w_{\rm reg}$,
\begin{eqnarray}
W\left[w_\sigma,\,w_{\rm reg}\right](\eta)=-\left(\sum_{k=0}^\infty a_{k,\sigma}\,\eta^{k+\sigma}\right)\left(\sum_{j=0}^\infty jb_j\,(1-\eta)^{j-1}\right) & &   \nonumber  \\
& & \hspace{-200pt}-\,\left(\sum_{k=0}^\infty (k+\sigma) a_{k,\sigma}\,\eta^{k-1+\sigma}\right)\left(\sum_{j=0}^\infty b_j\,(1-\eta)^j\right),  \label{iii6}
\end{eqnarray}
at an arbitrarily chosen point of the interval $-1<\eta<1$. From the computational point of view, an interesting choice of $\eta$ seems to be $\eta=1/2$, in which case
\begin{equation}
W\left[w_\sigma,\,w_{\rm reg}\right](1/2)=-\sum_{l=0}^\infty 2^{-l}\,(l+1)\left(\sum_{j=0}^{l+1-\sigma}a_{l+1-\sigma-j, \sigma}\,b_j\right)\,.  \label{iii7}
\end{equation}
The set of coefficients $\{a_{k,\sigma}\}$ and $\{b_j\}$  are solutions  of the difference equations (\ref{ii6}) and (\ref{ii10}), respectively. According to the Perron-Kreuser theorem on difference equations \cite{kreu,per1,per2},
\begin{equation}
\limsup_{k\to\infty}\left(|a_{k,\sigma}|\right)^{1/k}=1\,,   \qquad
\limsup_{j\to\infty}\left(|b_j|\right)^{1/j}=2^{-1}\,,   \label{iii8}
\end{equation}
that is, for any given $\varepsilon>0$, constants $C_a$ and $C_b$ can be found such that
\begin{equation}
|a_{k,\sigma}|<C_a(1+\varepsilon)^k\,\qquad   |b_j|<C_b(2^{-1}+\varepsilon)^j \qquad {\rm for\; any}\quad k,j\geq 0\,. \label{iii9}
\end{equation}
This makes evident that the series in the right hand side of Eq. (\ref{iii7}) converges as fast as the geometric series $\sum_{l=0}^\infty 2^{-l}$.
We consider, however, that a better choice of the value of $\eta$ in Eq. (\ref{iii6}) is $\eta=0$. In this case,
\begin{equation}
W\left[w_\sigma,\,w_{\rm reg}\right](0)=-\,\delta_{\sigma,0}\,\left(\sum_{j=0}^\infty j\,b_j\right) - \delta_{\sigma,1}\,\left(\sum_{j=0}^\infty b_j\right)\,.  \label{iii10}
\end{equation}
Obviously, the cancelation of this Wronskian occurs when either $dw_{\rm reg}/d\eta$ or $w_{\rm reg}$ vanish at the origin, as it occurs for respectively even or odd functions of $\eta$. Needless to say, the right hand side of (\ref{iii10}) depends on the variable $z$, introduced in (\ref{ii4}), through the coefficients $b_j$. Therefore, the problem of finding the eigenvalues of the angular spheroidal equations, for given $c^2$ and $m$, reduces to the determination of the zeros of the function
\begin{equation}
\mathcal{W}_\sigma(z)\equiv -\,W\left[w_\sigma,\,w_{\rm reg}\right](0)=\delta_{\sigma,0}\,\left(\sum_{j=0}^\infty j\,b_j(z)\right) + \delta_{\sigma,1}\,\left(\sum_{j=0}^\infty b_j(z)\right)\,,  \label{iii11}
\end{equation}
where we have indicated the dependence of the $b_j$ on $z$.

Different procedures can be used in the determination of the zeros of $\mathcal{W}_\sigma(z)$. A useful iterative method is the  Newton one. Starting with an initial approximate value, $z^{(0)}$, of a certain zero, repeated application of the algorithm
\begin{equation}
z^{(i+1)}=z^{(i)}-\frac{\mathcal{W}_\sigma(z^{(i)})}{\mathcal{W}_\sigma^{\prime}(z^{(i)})}   \label{iii12}
\end{equation}
allows one to get the value of the zero with the desired accuracy. The extended Newton method, which uses
\begin{equation}
z^{(i+1)}=z^{(i)}-\frac{\mathcal{W}_\sigma^\prime(z^{(i)})\pm \left[\left(\mathcal{W}_\sigma^\prime(z^{(i)})\right)^2-2\mathcal{W}_\sigma(z^{(i)})\,\mathcal{W}_\sigma^{\prime\prime}(z^{(i)})\right]^{1/2}} {\mathcal{W}_\sigma^{\prime\prime}(z^{(i)})}\,,   \label{iii13}
\end{equation}
is even more efficient. For the first and second derivatives of $\mathcal{W}_\sigma(z)$ with respect to $z$ we have the expressions
\begin{eqnarray}
\mathcal{W}_\sigma^\prime(z)&=& \delta_{\sigma,0}\,\left(\sum_{j=0}^\infty j\,b_j^\prime(z)\right) + \delta_{\sigma,1}\,\left(\sum_{j=0}^\infty b_j^\prime(z)\right)\,,  \label{iii14}  \\
\mathcal{W}_\sigma^{\prime\prime}(z)&=& \delta_{\sigma,0}\,\left(\sum_{j=0}^\infty j\,b_j^{\prime\prime}(z)\right) + \delta_{\sigma,1}\,\left(\sum_{j=0}^\infty b_j^{\prime\prime}(z)\right)\,.  \label{iii15}
\end{eqnarray}
The first and second derivatives of the coefficients $b_j(z)$ are easily obtained by means of the recurrence relations
\begin{eqnarray}
b_0^\prime(z)=0\,,\quad 2j(j+m)\,b_j^\prime(z)&=&\left[(j-1)(j+2m)-z+c^2\right]b_{j-1}^\prime(z) \nonumber   \\
& & \hspace{-40pt} -\,2c^2\,b_{j-2}^\prime(z)+c^2\,b_{j-3}^\prime(z)-b_{j-1}(z)\,,  \label{iii16} \\
b_0^{\prime\prime}(z)=0\,,\quad 2j(j+m)\,b_j^{\prime\prime}(z)&=&\left[(j-1)(j+2m)-z+c^2\right]b_{j-1}^{\prime\prime}(z)  \nonumber \\
& & \hspace{-40pt} -\,2c^2\,b_{j-2}^{\prime\prime}(z)+c^2\,b_{j-3}^{\prime\prime}(z)-2b_{j-1}^\prime(z)\,,   \label{iii17}
\end{eqnarray}
stemming from (\ref{ii10}). From these difference equations, inequalities analogous to the second one in (\ref{iii9}) can be deduced. Such inequalities guarantee the convergence, as fast as the geometric series $\sum_{j=0}^\infty 2^{-j}$, of the series in the right hand sides of Eqs. (\ref{iii11}), (\ref{iii14}), and (\ref{iii15}).

We have applied the procedure just described for obtaining the behaviour of the lowest eigenvalues of the spheroidal equation when the parameter $c^2$ varies in the interval $[-10, 10]$ and for values of $m=$0, 1, and 2. The results are shown in Figs. 1 to 3. (Remember that $\lambda_{m,n}=z_{m,n}+m(m+1)$.) A glance at Fig. 1 suggests a quasi-confluence of  trajectories of the eigenvalues $\lambda_{0,0}$ and $\lambda_{0,1}$ for sufficiently large negative values of $c^2$. One may conjecture that, for larger negative values of $c^2$ other pairs of trajectories, those of $\lambda_{0,2j}$ and $\lambda_{0,2j+1}$, present such quasi-confluence. Table 1 shows that this is the case, and that a similar phenomenon occurs for other values of $m$.
\begin{figure}
\begin{center}
\vspace{1cm}
\resizebox{10cm}{!}{\includegraphics{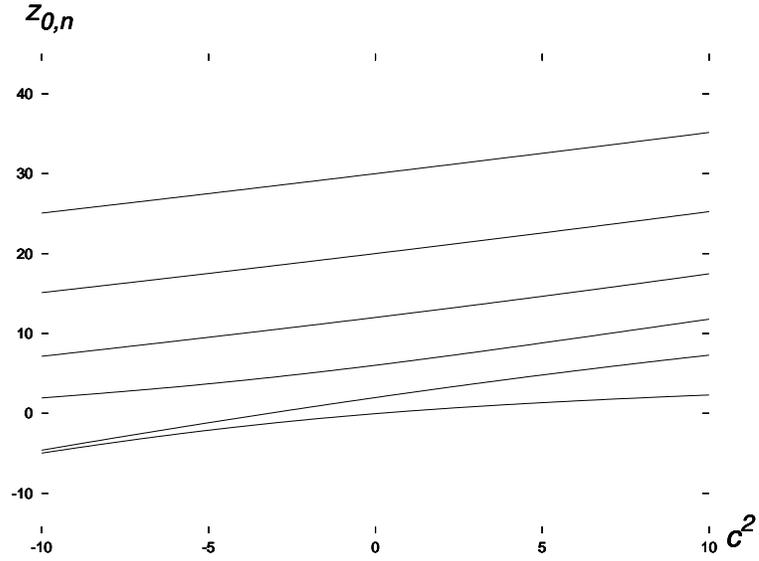}}
\end{center}
\caption{Trajectories of the lowest eigenvalues of the spheroidal equation with $m=0$ as $c^2$ varies in the interval $[-10, 10]$.} \label{m0}
\end{figure}
\begin{figure}
\begin{center}
\vspace{1cm}
\resizebox{10cm}{!}{\includegraphics{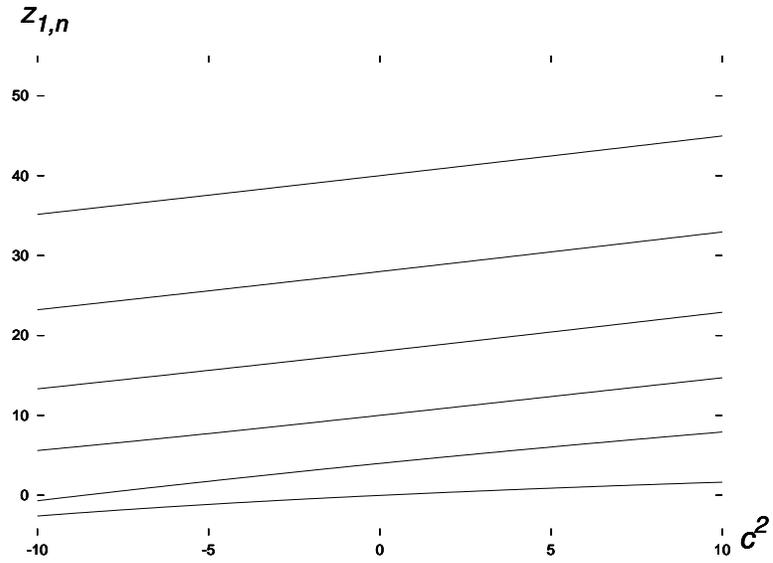}}
\end{center}
\caption{The same as in Figure 1, for $m=1$. Notice that, according to Eq. (\ref{ii4}), $\lambda_{1,n}=z_{1,n}+2$.} \label{m1}
\end{figure}
\begin{figure}
\begin{center}
\vspace{1cm}
\resizebox{10cm}{!}{\includegraphics{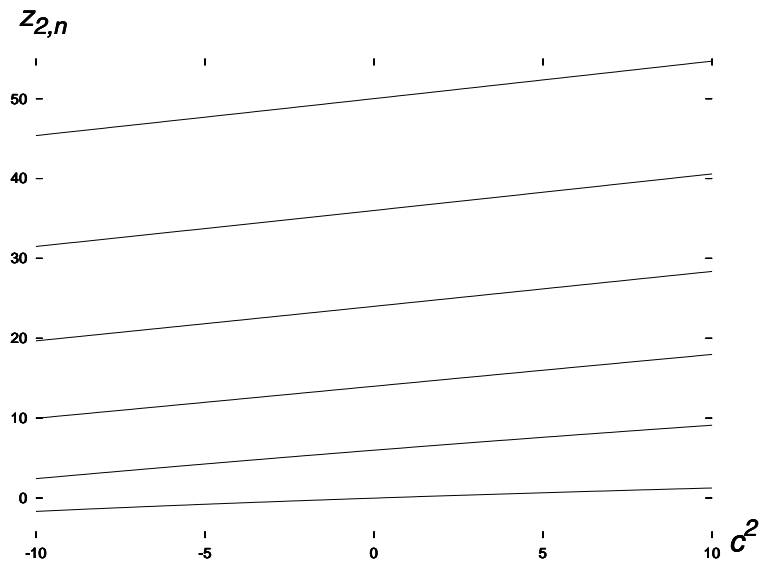}}
\end{center}
\caption{The same as in Figure 1, for $m=2$. In view of Eq. (\ref{ii4}), $\lambda_{2,n}=z_{2,n}+6$.} \label{m2}
\end{figure}
\begin {table}
\caption {Lowest eigenvalues of the (oblate) spheroidal equation for several values of $m$ and $c^2$.}
\begin{center}
\begin {tabular}{|l|r|r|r|}
\hline
\ & $m=0,\; c^2=-100$ & $m=1,\; c^2=-200$ & $m=2,\; c^2=-300$ \\
\hline
$\lambda_{m,m+5}$ & $-15.328144254756$ & $-51.05126046795$ & $-83.77105335717$ \\
$\lambda_{m,m+4}$ & $-16.065564650326$ & $-51.08618015853$ & $-83.77516906231$ \\
$\lambda_{m,m+3}$ & $-45.483938701812$ & $-95.57183718390$ & $-138.78472876574$ \\
$\lambda_{m,m+2}$ & $-45.489793371378$ & $-95.57199196249$ & $-138.78474405855$ \\
$\lambda_{m,m+1}$ & $-81.027938023746$ & $-145.51102178558$ & $-199.22477209684$ \\
$\lambda_{m,m}$   & $-81.027943944958$ & $-145.51102194107$ & $-199.22477211250$ \\
\hline
\end {tabular}
\end{center}
\end {table}

In order to compare with results published by other authors, we have applied our method to the computation of $\lambda_{m,n}$ for a sample of values of the parameters considered by Li {\em et al.} \cite{lilw} and by Ogburn {\em et al.} \cite{ogbu}. The comparison, shown in Table 2, allows one to conclude that, for moderate real values of $c^2$,  the procedure used in Ref. \cite{lilw} is more reliable than the  finite difference algorithm of Ref. \cite{ogbu}.
\begin {table}
\caption {Comparison of results obtained by using different procedures in the determination of the eigenvalues $\lambda_{m,n}$ of the spheroidal equation with real $c^2$.}
\begin{center}
\begin {tabular}{c c c l l l }
\hline
 $c^2$ & $m$ & $n$ & Ref. \cite{lilw} & Ref. \cite{ogbu} & this work\\
\hline
$-$1 & 4 & 11 & 131.5600809 & 131.560080918303  & 131.56008091940694  \\
0.1 & 2 & 2 & 6.014266314 & 6.014266356124070  & 6.0142663139415926  \\
1 & 1 & 1   & 2.195548355 & 2.195612369653500  & 2.1955483554130039  \\
1 & 2 & 2   & 6.140948992 & 6.140948969717170  & 6.1409489918576905  \\
1 & 2 & 5   & 30.43614539 & 30.436145317468500 & 30.436145388713659  \\
4 & 1 & 1   & 2.734111026 & 2.73415086499219   & 2.7341110256122556  \\
4 & 2 & 2   & 6.542495274 & 6.54249530312951   & 6.5424952743905705  \\
16 & 1 & 1  & 4.399593067 & 4.399599760664940  & 4.3995930671655061  \\
16 & 2 & 5  & 36.99626750 & 36.996267483327900 & 36.996267500847930  \\
\hline
\end {tabular}
\end{center}
\end {table}

Obviously, the procedure is also applicable in the case of complex $c^2$. Table 3 shows our results for different values of $c$, $m$, and $n$ considered in Ref \cite{lilw}. As it can be seen, the eigenvalues given by Li {\em et al.} are confirmed. Nevertheless, in the neighbourhood of each one of those eigenvalues, we have found another one, reported also in Table 3. This result is not surprising, because the values of $c$ considered are the approximations found by  Oguchi \cite{oguc} to what he calls ``the branch points of the eigenvalues as functions of $c$", that is, in our formalism, values of $c$ for which a double zero of $\mathcal{W}_\sigma(z)$ exists. According to the results of Skorokhodov and Khristoforov \cite{skor}, there is a double eigenvalue $\lambda=1.705180091+4.220186348\,i$ for $c=1.824770749+2.601670693\,i$. This is a much better approximation to the branch point unveiled by Oguchi. As an illustration of what happens in the vicinity of those values of $c$, we present in Figure 4 a modulus-phase plot of the function $\mathcal{W}_0(z)$ for $c=1.824770+2.601670\,i$ and $m=0$, the first of the cases considered in Table 3. The two eigenvalues reported in the table appear as zeros of $\mathcal{W}_0(z)$. As the value of $c$ moves from the approximation to the branch point found by Oguchi towards the more precise value given by Skorokhodov and Khristoforov, the two zeros of $\mathcal{W}_0(z)$ shown in Fig.~4 approach to each other and eventually collide at a point in the close neighbourhood of the saddle point of $\mathcal{W}_0(z)$ suggested by its modulus-phase plot. Similar plots of $\mathcal{W}_\sigma(z)$ are obtained for the other cases in Table 3.

\begin {table}
\caption {Pairs of eigenvalues $\lambda_{m,n}$ of the spheroidal equation for the complex values of $c$ given in Ref \cite{oguc} as corresponding to ``branch points".}
\hspace{-20pt}
\begin {tabular}{c c c c c }
\hline
 $c$ & $m$ & $n$ & Ref. \cite{lilw} & this work \\
\hline
$1.824770+2.601670\,i$ & 0 & 0 & $1.701836+4.219998\,i$ & $1.701836497+4.219997758\,i$ \\
 & & 2 & & $1.708523909+4.220369152\,i$ \\
 & & & & \\
$2.094267+5.807965\,i$ & 0 & 0 & $1.993901+8.576325\,i$ & $1.993900944+8.576324731\,i$ \\
 & & 4 & & $2.003141811+8.581103855\,i$ \\
 & & & & \\
$5.217093+3.081362\,i$ & 0 & 2 & $23.91023+18.74194\,i$ & $23.91033400+18.74184255\,i$ \\
 & & 4 &  & $23.92132979+18.74479980\,i$ \\
 & & & & \\
$3.563644+2.887165\,i$ & 0 & 1 & $10.13705+11.12216\,i$ & $10.13704735+11.12217988\,i$ \\
 & & 3 & & $10.14462729+11.12098765\,i$ \\
 & & & & \\
$1.998555+4.097453\,i$ & 1 & 1 & $2.919098+6.134851\,i$ & $2.919095372+6.134851876\,i$ \\
 & & 3 & & $2.911544002+6.133045176\,i$ \\
 & & & & \\
$3.862833+4.492300\,i$ & 1 & 2 & $12.19691+16.24534\,i$ & $12.19691647+16.24534182\,i$ \\
 & & 4 & & $12.20527134+16.24281200\,i$ \\
 & & & & \\
$2.136987+5.449457\,i$ & 2 & 2 & $6.098946+7.684379\,i$ & $6.098961456+7.684332819\,i$ \\
 & & 4 & & $6.106119819+7.685191032\,i$ \\
\hline
\end {tabular}
\end {table}
\begin{figure}
\begin{center}
\vspace{1cm}
\resizebox{11cm}{!}{\includegraphics{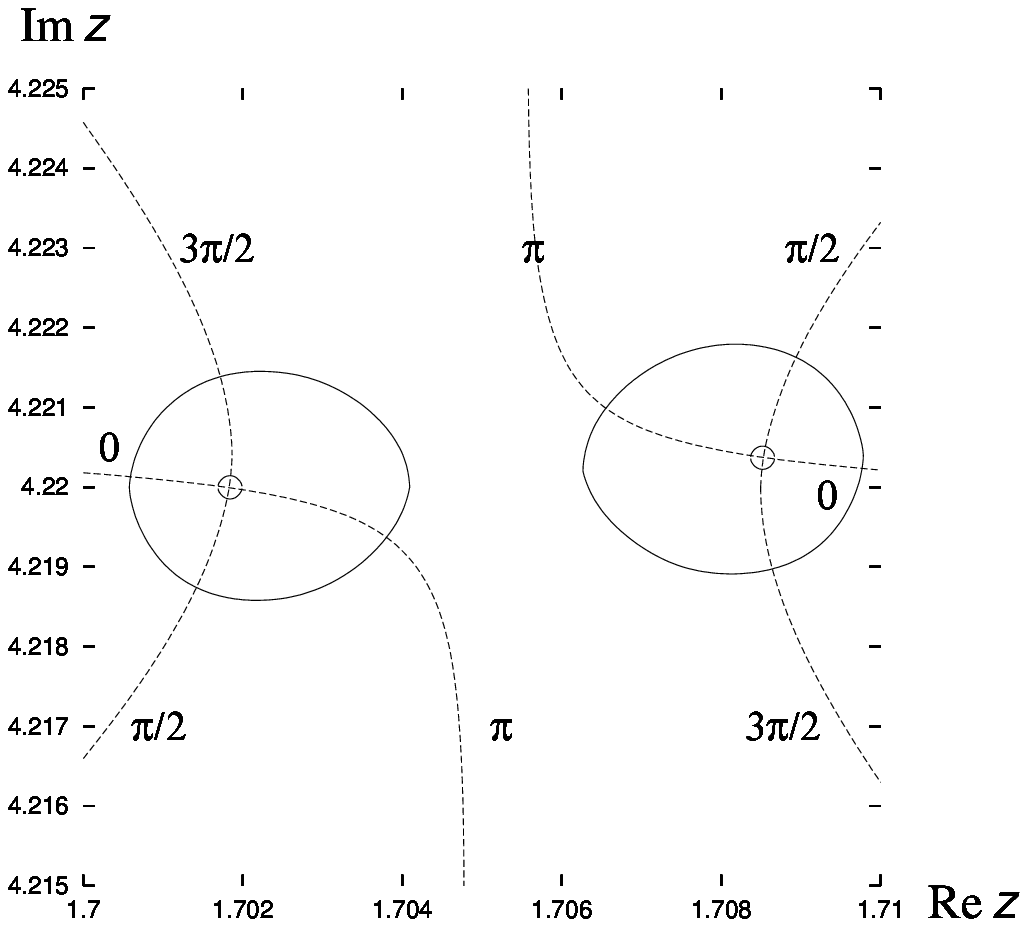}}
\end{center}
\caption{Modulus-phase plot of $\mathcal{W}_0(z)$ for $c=1.824770+2.601670\,i$, in the neighbourhood of a ``branch point", and $m=0$. Continuous and dashed lines are used to represent, respectively, the constant-modulus and constant-phase loci. Only the constant-modulus lines corresponding to $|\mathcal{W}_0(z)|=10^{-6}$ and $10^{-7}$  and the constant-phase lines for $\arg\mathcal{W}_0(z)=0,\; \pi/2,\; \pi$ and $3\pi/2$ have been drawn. }
\label{modfas}
\end{figure}

A comment concerning the values of the label $n$ of $\lambda_{m,n}$ reported in Table 3 is in order. For real or pure imaginary $c$, i.e. for real $c^2$, the eigenvalues $\lambda_{m,n}$ for given $m$ are real and can be ordered by increasing value. The label $n$ reflects that order. For complex $c^2$, instead, the values of $\lambda(c)$ become complex and such ordination is no more possible. Nevertheless, a label $n$ can be assigned to those complex values of $\lambda$, as done by Skorokhodov and Khristoforov. By keeping constant the real part of $c$ and continuously decreasing its imaginary part, $\lambda(c)$ describes, in the complex $\lambda$-plane, a trajectory which intersects the real $\lambda$-axis at a certain $\lambda_{m,n}$ for $\Im c=0$. This label $n$ can be attached to the whole trajectory described by $\lambda(c)$ as $c$ varies in the complex plane. The paper by Skorokhodov and Khristoforov contains a very lucid discussion of those trajectories and shows that the branch points $c_s$ correspond to singular values of $c$ such that $\lambda_{m,n}(c_s)=\lambda_{m,n+2p}(c_s)$, with $p=1, 2, \ldots $.

\section{The eigenfunctions}

Once the eigenvalues $\lambda_{m,n}$ have been calculated, the corresponding eigenfunctions, in the interval $0\leq\eta\leq 1$, can be obtained immediately by means of the series expansion
\begin{equation}
S_{m,n}(c,\eta)=\mathcal{N}\,e^{i\theta}\,(1-\eta^2)^{m/2}\,\sum_{j=0}^\infty b_j\,(1-\eta)^j\,,   \label{iv1}
\end{equation}
the coefficients $b_j$ being given by the recurrence relation (\ref{ii10}) with $z=\lambda_{m,n}-m(m+1)$. The normalization constant $\mathcal{N}$ should be adjusted to the normalization scheme preferred. A discussion of  the different normalizations used in the literature can be found in the paper by Kirby \cite{kirb}, where the advantage of the unit normalization
\begin{equation}
\int_{-1}^1 |S_{m,n}(c,\eta)|^2\,d\eta=1   \label{iv2}
\end{equation}
is made evident. By choosing this normalization, one has
\begin{equation}
\mathcal{N}=\left[ \sum_{k=0}^m (-1)^k\,2^{m+1-k}\,{m \choose k}\,\sum_{l=0}^\infty \frac{\sum_{j=0}^l b_j\,b_{l-j}^*}{l+m+k+1}\right]^{-1/2}\,,  \label{iv3}
\end{equation}
where the asterisk indicates complex conjugation.
The constant phase $\theta$ in the right hand side of (\ref{iv1}) may be taken at will. In the case of real $c^2$, it is natural to take $\theta =0$. For complex $c^2$, $\theta$ can be chosen in such a way that $S_{m,n}$ becomes real at $\eta=0$, or at $\eta=1$, or at any other point. Needless to say, $S_{m,n}(c,-\eta)=\pm S_{m,n}(c,\eta)$, according to the even or odd nature of $S_{m,n}$.

Figures 5 and 6 show two examples of the application of the method to the computation of spheroidal angular wave functions in the case of real $c^2$. The first one is an even prolate angular wave function of parameters $c=3$, $m=0$ and $n=2$, and eigenvalue $\lambda_{0,2}=11.192938649526784$, a case considered in Table III of Ref. \cite{lilw} (with the Flammer \cite{flam} normalization scheme). The second one is an odd oblate angular wave function corresponding to the first of the cases considered in our Table 2. In both figures, the functions have been normalized according to Eq. (\ref{iv2}). We have applied our procedure also in some cases of complex $c$. Figure 7 shows the real and imaginary parts and the squared modulus of the wave function in a case considered by Falloon {\em et al.} \cite{fall}, namely the first one in their Table 2. The parameters are $c=1\! +\! i$, $m=0$ and $n=0$, and the eigenvalue, in our notation, is $\lambda_{0,0}=0.059472769735031+0.662825122194600\,i$. (We give here the eigenvalue with only 15 decimal digits, but our procedure is able to reproduce the 25 decimal digits given in Ref. \cite{fall} and to obtain even more.) Finally, Figures 8, 9 and 10 correspond to the second of the cases considered in Table 2 of Ref. \cite{ogbu}, of parameters $c=20(1\! +\! i)$, $m=0$, $n=3$, and eigenvalue $\lambda_{0,3}=58.226714354344554 + 60.025615481720256\,i$. (Notice the discrepancy, in the six last digits of both real and imaginary parts, with the value given in Ref. \cite{ogbu}.)
\begin{figure}
\begin{center}
\vspace{1cm}
\resizebox{10cm}{!}{\includegraphics{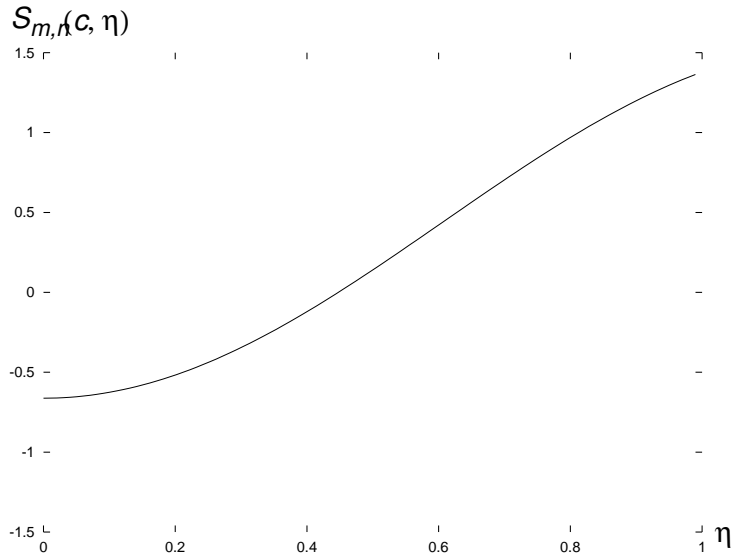}}
\end{center}
\caption{Prolate spheroidal angular wave function of parameters $c=3$, $m=0$ and $n=2$, corresponding to the eigenvalue  $\lambda_{0,2}=11.192938649526788$. Since $S_{0,2}(\eta)$ is an even function, we have omitted its representation in the interval $-1\leq\eta<0$. The normalization adopted is that prescribed in Eq. (\ref{iv2}).}
\label{wfeven}
\end{figure}
\begin{figure}
\begin{center}
\vspace{1cm}
\resizebox{10cm}{!}{\includegraphics{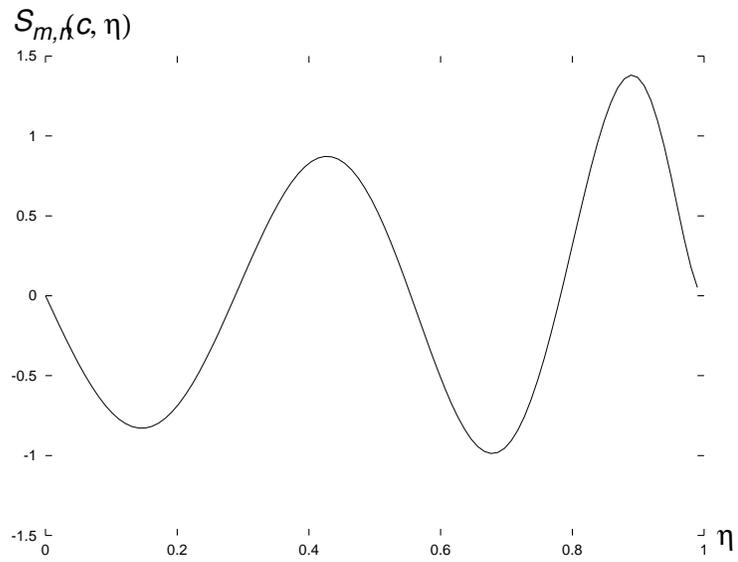}}
\end{center}
\caption{Oblate spheroidal angular wave function of parameters $c=i$, $m=4$ and $n=11$, with eigenvalue  $\lambda_{4,11}=131.56008091940694$. The function is an odd one. It has been normalized as in Eq. (\ref{iv2}).}
\label{wfodd}
\end{figure}
\begin{figure}
\begin{center}
\vspace{1cm}
\resizebox{10cm}{!}{\includegraphics{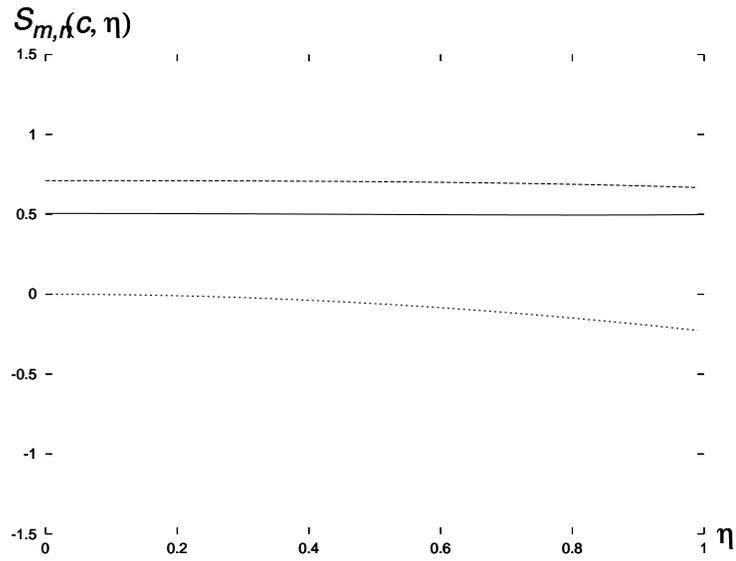}}
\end{center}
\caption{Real and imaginary parts and squared modulus of the angular spheroidal wave function of parameters  $c=1\! +\! i$, $m=0$ and $n=0$, and eigenvalue $\lambda_{0,0}=0.059472769735031+0.662825122194600\, i$, normalized to unit, as in Eq. (\ref{iv2}). Dashed and dotted lines are used to represent, respectively, the real and imaginary parts of $S_{0,0}(1\! +\! i, \eta)$, and a solid line for its squared modulus. The arbitrary phase $\theta$ in the right hand side of Eq. (\ref{iv1}) has been fixed in such a way that the wave function becomes real at the origin. Only the interval $0\leq\eta\leq 1$ has been considered. Needless to say, $S_{0,0}(c, -\eta)=S_{0,0}(c, \eta)$.}
\label{fallon}
\end{figure}
\begin{figure}
\begin{center}
\vspace{1cm}
\resizebox{10cm}{!}{\includegraphics{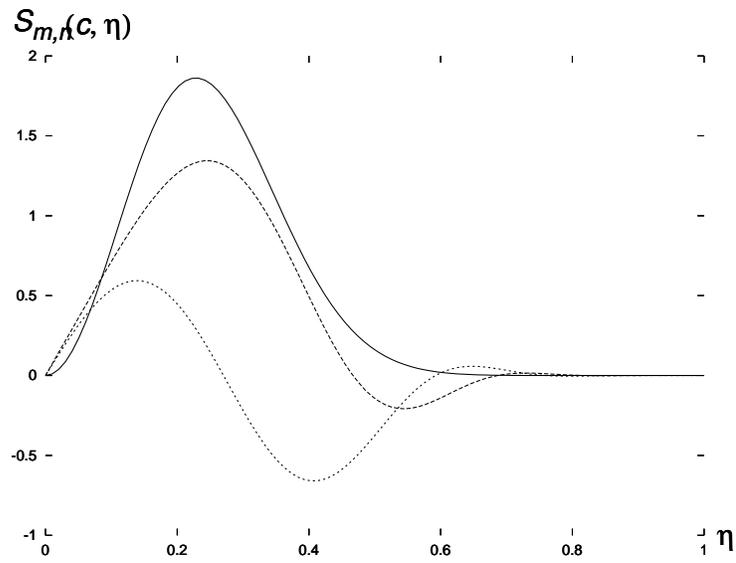}}
\end{center}
\caption{Real and imaginary parts and squared modulus of the angular spheroidal wave function of parameters  $c=20(1\! +\! i)$, $m=0$ and $n=3$, and eigenvalue $\lambda_{0,3}=58.226714354344554 + 60.025615481720256\, i$. The meaning of the lines and the normalization is the same as in Fig. 7. For the arbitrary phase $\theta$ in Eq. (\ref{iv1}) we have chosen the value $\theta=0$. Of course, $S_{0,3}(c, -\eta)=-\,S_{0,3}(c, \eta).$}
\label{barrowes0}
\end{figure}
\begin{figure}
\begin{center}
\vspace{1cm}
\resizebox{10cm}{!}{\includegraphics{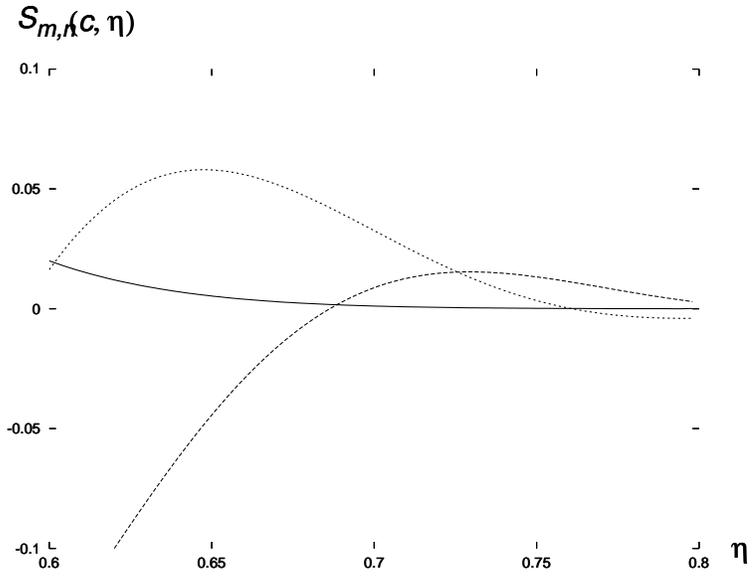}}
\end{center}
\caption{Magnification of the interval $0.6\leq\eta\leq 0.8$ of Fig. 8.}
\label{barrowes1}
\end{figure}
\begin{figure}
\begin{center}
\vspace{1cm}
\resizebox{10cm}{!}{\includegraphics{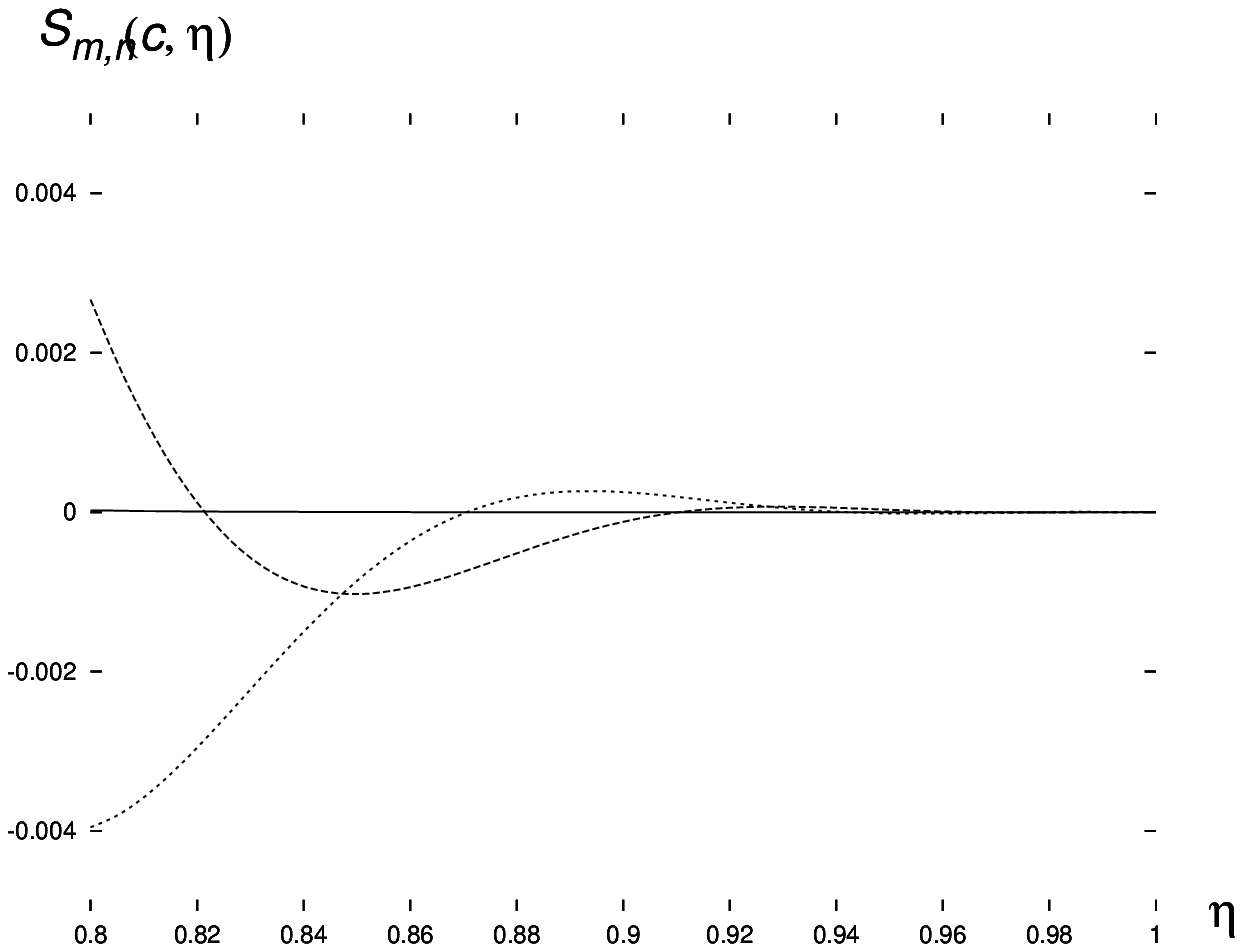}}
\end{center}
\caption{Magnification of the interval $0.8\leq\eta\leq 1$ of Fig. 8.}
\label{barrowes2}
\end{figure}

\section{Conclusions}

We have developed a rarely used method which allows to find the eigenvalues and eigenfunctions of the angular spheroidal equation. Instead of having recourse to comparison with other differential equations, the procedure deals with a direct solution, expressed in the form of a convergent series. Requiring it to be regular gives the eigenvalues, which appear as the zeros of a one variable function, $\mathcal{W}_\sigma(z)$. This function and its derivatives $\mathcal{W}_\sigma^\prime(z)$ and $\mathcal{W}_\sigma^{\prime\prime}(z)$ with respect to the variable $z$ can be computed, to the desired precision, by summing rapidly convergent series. This fact makes possible the application of the extended Newton method for the determination of the zeros of $\mathcal{W}_\sigma(z)$, i. e., the eigenvalues of the spheroidal equation. Then, the computation of the corresponding eigenfunctions, conveniently normalized, becomes trivial. For the normalization, one benefits from the fact that the squared modulus of the wave function can be integrated algebraically.

The fact that, for given $c$ and $m$, the eigenvalues are the zeros of an easily computable function, $\mathcal{W}_\sigma(z)$, makes possible to get an initial approximate location of all of them by a tabulation or a graphical representation of that function. Repeated application of the extended Newton method allows then to calculate the eigenvalues with the desired accuracy.

We have shown the applicability of the method not only in the cases of prolate (real $c$) and oblate (imaginary $c$) spheroidal wave equations, but also when $c$ is complex. The procedure provides in all cases a very precise determination of the eigenvalues. This has allowed us to resolve the quasi-confluence of pairs of even-odd eigenvalues for large imaginary values of $c$ (Table 1) and of pairs of even-even or odd-odd eigenvalues for complex values of $c$ in the neighbourhood of ``branch points" (Table 3).

The efficiency of the procedure proposed in this paper is subordinate to the capability of computing $\mathcal{W}_\sigma(z)$ with sufficient accuracy. The convergence of the series in the right hand side of Eq. (\ref{iii11}) is guaranteed in all cases, since $|b_n|\sim 2^{-n}$ for all $n$ larger than a certain $N$, and the series may replaced by a sum up to say $j=j_{\rm max}$. Nevertheless, for large values of $c$ and/or $\lambda_{m,n}$, the coefficients $b_j$ increase (in modulus) rapidly with $j$ before starting to decrease, and the adequate value of $j_{\rm max}$ may become very large.
Even worse, the values of the terms to be summed may cover so many orders of magnitude that the resulting sum is not reliable, unless many significant digits are carried along the computation. This drawback is not outside other procedures. However, the algorithms proposed by Kirby \cite{kirb} and by Ogburn {\em et al.} \cite{ogbu}, and procedures collected in Refs.~\cite{osi2}, seem to be able to tackle the issue properly. Asymptotic methods \cite{barr,rokh} have also been forwarded for the mentioned cases of large values of $c$ and/or $\lambda_{m,n}$.

\section*{Acknowledgements}

The work has been supported by Departamento de Ciencia, Tecnolog\'{\i}a y Universidad del Gobierno de Arag\'on (Project 226223/1) and Ministerio de Ciencia
e Innovaci\'on (Project MTM2015-64166)

\bibliographystyle{amsplain}

\end{document}